\documentclass[12pt]{article}

\begin{document}

\begin{titlepage}

\title{\bf Complex Lagrangian mechanics \\ with constraints}
\author{ Mehmet Tekkoyun \footnote{tekkoyun@pau.edu.tr} \\
 {\small Department of Mathematics, Pamukkale University,}\\
{\small 20070 Denizli, Turkey}\\
Ali G\"{o}rg\"{u}l\"{u} \footnote{agorgulu@ogu.edu.tr} \\
 {\small Department of Mathematics, Eski\c{s}ehir Osmangazi University,}\\
 {\small 26480  Eski\c{s}ehir, Turkey}}
\date{\today}
\maketitle
\begin{abstract}

In this study, it is generalized the concept of Lagrangian
mechanics with constraints to complex case. To be beginning , it
is considered a K\"{a}hlerian manifold as a velocity-phase space.
Then a non-holonomic constraint is given by 1-form on it. If the
form is closed, it is found that the constraint is (locally)
holonomic. In the result, complex analogous of some topics in
constrained Lagrangian mechanical system is concluded.


\end{abstract}

\section{Introduction}

Modern differential geometry provides a fundamental framework for studying
Lagrangian mechanics. In recent years, there are many studies as some
articles in \cite{crampin,tekkoyun1,tekkoyun2} and books in \cite
{deleon1,deleon2} about differential geometric methods in mechanics. It is
well known that the dynamics of Lagrangian formalisms is characterized by a
suitable vector field defined on the tangent bundles which are phase-spaces
of velocities of a given configuration manifold. If $Q$ is an $m$%
-dimensional configuration manifold and $L:TQ\rightarrow \mathbf{R}$ is a
regular Lagrangian function, then there is a unique vector field $\xi _{L}$
on $TQ$ such that dynamical equations
\begin{equation}
i_{\xi _{L}}\Phi _{L}=dE_{L},  \label{1.1}
\end{equation}
where $\Phi _{L}$ is the symplectic form and $E_{L}$ is the energy
associated to $L$. The Euler-Lagrange vector field $\xi _{L}$ is a semispray
(or \textit{second order differential equation}) on $Q$ since its integral
curves are the solutions of the Euler-Lagrange equations given by
\begin{equation}
\frac{d}{dt}\frac{\partial L}{\partial \dot{q}^{i}}-\frac{\partial L}{%
\partial q^{i}}=0,  \label{1.2}
\end{equation}
where $q^{i}$ and $(q^{i},{\dot{q}}^{i}),1\leq i\leq m,$ are coordinate
system of $Q$ and $TQ.$ The triple, either $(TQ,\Phi _{L},\xi _{L})$ or $%
(TQ,\Phi _{L},E_{L})$ is called \textit{Lagrangian mechanical system} on the
tangent bundle $TQ.$ Assume that $(TQ,\Phi _{L})$ is symplectic manifold and
$\overline{\omega }=\{\omega _{1},...,\omega _{r}\}$ is a system of
constraints on $TQ.$ We call to be a \textit{constraint} on $TQ$\ to a
non-zero 1-form $\omega =\wedge ^{a}\omega _{a}$ on $TQ,$ such that $\wedge
^{a}$ are Lagrange multipliers$.$ We call $(TQ,\Phi _{L},E_{L},\overline{%
\omega })$ a \textit{regular Lagrangian system with constraints}. The
constraints $\overline{\omega }$ are said to be classical constraints if the
1-forms $\omega _{a},1\leq a\leq r,$ are basic. Then holonomic classical
constraints define foliations on the configuration manifold $Q,$ but
holonomic constraints also admit foliations on the phase space of velocities
$TQ.$ As real studies, generally a curve $\alpha $ satisfying the Euler
Lagrange equations for Lagrangian energy $E_{L}$ will not satisfy the
constraints. It must be that some additional forces (or \textit{canonical
constraint forces}) act on the system in addition to the \textit{force} $%
dE_{L}$ for a curve $\alpha $ to satisfy the constraints$.$ It is said that
the quartet $(TQ,\Phi _{L},E_{L},\overline{\omega })$ defines a \textit{%
mechanical system with constraints} if the vector field $\xi $ given by the
equations of motion
\begin{equation}
i_{\xi }\Phi _{L}=dE_{L}+\wedge ^{a}\omega _{a},\,\,\,\,\,\,\,\,\,\omega
_{a}(\xi )=0,  \label{1.3}
\end{equation}
is a semispray. Then, it is given Euler-Lagrange equations with constraints
as follows:
\begin{equation}
\frac{\partial L}{\partial q^{i}}-\frac{d}{dt}\frac{\partial L}{\partial
\stackrel{.}{q}^{i}}=\wedge ^{a}(\omega _{a})_{i},  \label{1.4}
\end{equation}

The purpose of this study is to make a contribution to the modern
development of Lagrangian formalisms of classical mechanics in terms of
differential-geometric methods on differentiable manifolds. So, we obtain
complex Euler-Lagrange equations with constraints on the K\"{a}hlerian
manifold. In the conclusion section, geometrical and mechanical results of
constrained mechanical system have been given. The first of them is that if
the distribution on $TM$ is integrable, a system of constraints is
holonomic. The second is that constrained Lagrange energy is conserved.

The present paper is structured as follows. In sections 2, it is recalled
complex and K\"{a}hlerian manifolds, and also Euler-Lagrange equations on
K\"{a}hlerian manifolds. In section 3, complex Euler-Lagrange equations with
constraints on K\"{a}hlerian manifold are deduced. In the final section, the
geometrical and mechanical meaning of constrained complex mechanical system
was given.

\section{Preliminaries}

In this letter, all geometric objects are assumed to be differentiable and
the sum is taken over repeated indices$.$ Now then it is assumed $1\leq
i\leq m.$

\subsection{Complex manifolds}

Let $M$ be configuration manifold of real dimension $m.$ A tensor field $J$
on $TM$ is called an \textit{almost complex structure} on $TM$ if at every
point $p$ of $TM,$ $J$ is endomorphism of the tangent space $T_{p}(TM)$ such
that $J^{2}=-I.$ A manifold $TM$ with fixed almost complex structure $J$ is
called \textit{almost complex manifold}. Assume that $(x_{i})$ be
coordinates of $M$ and $(x_{i},\,y_{i})$ be a real coordinate system on a
neighborhood $U $ of any point $p$ of $TM.$ Also, let us to be $\{(\frac{%
\partial }{\partial x^{i}})_{p},(\frac{\partial }{\partial y^{i}})_{p}\}$
and $\{(dx^{i})_{p},(dy^{i})_{p}\}$ to natural bases over $\mathbf{R}$ of
tangent space $T_{p}(TM)$ and cotangent space $T_{p}^{*}(TM)$ of $TM$,
respectively$. $

Let $TM$ be an almost complex manifold with fixed almost complex structure $%
J.$ The manifold $TM$ is called \textit{complex manifold} if there exists an
open covering $\left\{ U\right\} $ of $TM$ satisfying the following
condition: There is a local coordinate system $(x_{i},\,y_{i})$ on each $U,$
such that
\begin{equation}
J(\frac{\partial }{\partial x_{i}})=\frac{\partial }{\partial y_{i}},\,J(%
\frac{\partial }{\partial y_{i}})=-\frac{\partial }{\partial x_{i}}.
\label{2.1}
\end{equation}

for each point of $U.$ Let $z_{i}=x_{i}+$i$\,y_{i},\,$i$=\sqrt{-1},$ be a
complex local coordinate system on a neighborhood $U$ of any point $p$ of $%
TM.$ We define the vector fields by
\begin{equation}
(\frac{\partial }{\partial z^{i}})_{p}=\frac{1}{2}\{(\frac{\partial }{%
\partial x^{i}})_{p}-\mathbf{i}(\frac{\partial }{\partial y^{i}})_{p}\},\,(%
\frac{\partial }{\partial \overline{z}^{i}})_{p}=\frac{1}{2}\{(\frac{%
\partial }{\partial x^{i}})_{p}+\mathbf{i}(\frac{\partial }{\partial y^{i}}%
)_{p}  \label{2.2}
\end{equation}

and the dual covector fields
\begin{equation}
\left( dz^{i}\right) _{p}=\left( dx^{i}\right) _{p}+\mathbf{i}%
(dy^{i})_{p},\,\left( d\overline{z}^{i}\right) _{p}=\left( dx^{i}\right)
_{p}-\mathbf{i}(dy^{i})_{p}  \label{2.3}
\end{equation}

which represent bases of the tangent space $T_{p}(TM)$ and cotangent space $%
T_{p}^{*}(TM)$ of $TM,$ respectively. Then the endomorphism $J$ is shown as
\begin{equation}
J(\frac{\partial }{\partial z_{i}})=\mathbf{i}\frac{\partial }{\partial z_{i}%
},\,J(\frac{\partial }{\partial \overline{z}_{i}})=-\mathbf{i}\frac{\partial
}{\partial \overline{z}_{i}}.  \label{2.4}
\end{equation}

The dual endomorphism $J^{*}$ of the cotangent space $T_{p}^{*}(TM)$ at any
point $p$ of manifold $TM$ satisfies $J^{*2}=-I,$ and is defined by
\begin{equation}
J^{*}(dz_{i})=\mathbf{i}dz_{i},\,J^{*}(d\overline{z}_{i})=-\mathbf{i}d%
\overline{z}_{i}.  \label{2.5}
\end{equation}

\subsection{Hermitian and K\"{a}hlerian manifolds}

A \textit{Hermitian metric} on an almost complex manifold with almost
complex structure $J$ is a Riemannian metric $g$ on $TM$ such that
\begin{equation}
g(JX,JY)=g(X,Y),  \label{2.6}
\end{equation}

for any vector fields $X$, $Y$ on $TM.$ An almost complex manifold $TM$ with
a Hermitian metric is called an \textit{almost Hermitian manifold}. If,
moreover, $TM$ is a complex manifold, then $TM$ is called a \textit{%
Hermitian manifold}.

Let further $TM$ be a 2m-dimensional real almost Hermitian manifold with
almost complex structure $J$ and Hermitian metric $g$. The triple $(TM,J,g)$
may be named an \textit{almost Hermitian structure}. We denote by $\chi (TM)$
the set of complex vector fields on $TM$ and by $\wedge ^{1}(TM)$ the set of
complex 1-forms on $TM.$ Let $(TM,J,g)$ be an almost Hermitian structure.
The 2-form defined by
\begin{equation}
\Phi (X,Y)=g(X,JY),\,\,\,\,\forall X,Y\in \chi (TM)  \label{2.7}
\end{equation}

is called the \textit{K\"{a}hlerian form} of $(TM,J,g).$

An almost Hermitian manifold is called \textit{almost K\"{a}hlerian} if its
K\"{a}hlerian form\textbf{\ }$\Phi $ is closed. If, moreover, $TM$ is
Hermitian, then $TM$ is called a K\"{a}hlerian manifold.

\subsection{Complex Euler-Lagrange Equations}

In this subsection, it is recalled complex Euler-Lagrange equations for
classical mechanics structured on K\"{a}hlerian manifold introduced in \cite
{tekkoyun1}.

Let $J$ be an almost complex structure on the K\"{a}hlerian manifold and $%
(z^{i},\overline{z}^{i})$ its complex coordinates. We call to be the
semispray to the vector field $\xi $ given by

\begin{equation}
\xi =\xi ^{i}\frac{\partial }{\partial z^{i}}+\overline{\xi }^{i}\frac{%
\partial }{\partial \overline{z}^{i}},\xi ^{i}=\stackrel{.}{z}^{i}=\overline{%
z}^{i},\overline{\xi }^{i}=\stackrel{.}{\xi }^{i}=\stackrel{..}{z}^{i}=%
\stackrel{.}{\overline{z}}^{i}.  \label{2.8}
\end{equation}
The vector field $V=J\xi $ is called \textit{Liouville vector field} on the
K\"{a}hlerian manifold. We call \textit{the kinetic energy} and \textit{the
potential energy of system} the maps given by $T,P:TM\rightarrow \mathbf{C}$
such that $T=\frac{1}{2}m_{i}(\overline{z}^{i})^{2}=\frac{1}{2}m_{i}(%
\stackrel{.}{z}^{i})^{2},P=m_{i}\mathbf{g}h,$ respectively, where $m_{i}$ is
mass of a mechanic system having $m$ particles, $\mathbf{g}$ is the gravity
acceleration and $h$ is the origin distance of the a mechanic system on the
K\"{a}hlerian manifold. Then it may be said to be \textit{Lagrangian function%
} the map $L:TM\rightarrow \mathbf{C}$ such that $L=T-P$ and also \textit{%
the energy function} associated $L$ the function given by $E_{L}=V(L)-L.$

The vertical derivation operator $i_{J}$ defined by

\begin{equation}
i_{J}\omega(Z_{1},Z_{2},...,Z_{r})=\textstyle\sum_{i=1}^r \omega
(Z_{1},...,JZ_{i},...,Z_{r}),  \label{2.9}
\end{equation}
where $\omega \in \wedge ^{r}TM,$ $Z_{i}\in \chi (TM).$ The exterior
differentiation $d_{J}$ is defined by

\begin{equation}
d_{J}=[i_{J},d]=i_{J}d-di_{J},  \label{2.10}
\end{equation}
where $d$ is the usual exterior derivation.

For almost complex structure $J$, the closed K\"{a}hlerian form is the
closed 2-form given by
\begin{equation}
\Phi _{L}=-dd_{J}L,  \label{2.11}
\end{equation}
such that
\[
d_{J}:\mathcal{F}(TM)\rightarrow \wedge ^{1}TM
\]

By means of (\ref{1.3}), \textit{complex Euler-Lagrange equations} on
K\"{a}hlerian manifold $TM$ is found the following as:

\begin{equation}
\mathbf{i}\frac{\partial }{\partial t}\left( \frac{\partial L}{\partial z^{i}%
}\right) -\frac{\partial L}{\partial z^{i}}=0,\,\,\,\,\,\,\,\,\,\,\,\,\,\,%
\mathbf{i}\frac{\partial }{\partial t}\left( \frac{\partial L}{\partial
\stackrel{.}{z}^{i}}\right) +\frac{\partial L}{\partial \stackrel{.}{z}^{i}}%
=0.  \label{2.12}
\end{equation}

\section{Complex Euler-Lagrange Equations with Constraints}

In this section, we shall obtain the version with constraints of complex
Euler-Lagrange equations for classical mechanics structured on K\"{a}hlerian
manifold introduced in \cite{tekkoyun1}.

Let $J$ be an almost complex structure on the K\"{a}hlerian manifold and $%
(z^{i},\overline{z}^{i})$ its complex coordinates.. Assume to be semispray
to the vector field $\xi $ given as:

\begin{equation}
\xi =\xi _{L}+\wedge ^{a}\omega _{a}=\xi ^{i}\frac{\partial }{\partial z^{i}}%
+\overline{\xi }^{i}\frac{\partial }{\partial \overline{z}^{i}}+\wedge
^{a}\omega _{a},\,\,1\leq a\leq r,  \label{3.1}
\end{equation}
The vector field determined by
\begin{equation}
V=J\xi _{L}=\mathbf{i}\xi ^{i}\frac{\partial }{\partial z^{i}}-\mathbf{i}%
\overline{\xi }^{i}\frac{\partial }{\partial \overline{z}^{i}},  \label{3.2}
\end{equation}
is called \textit{Liouville vector field }on the K\"{a}hlerian manifold $TM$%
. The closed 2-form given by $\Phi _{L}=-dd_{J}L$ such that
\begin{equation}
d_{J}=\mathbf{i}\frac{\partial }{\partial z^{i}}dz^{i}-\mathbf{i}\frac{%
\partial }{\partial \overline{z}^{i}}d\overline{z}^{i}:\mathcal{F}%
(TM)\rightarrow \wedge ^{1}TM.  \label{3.5}
\end{equation}
is found to be
\begin{eqnarray}
\Phi _{L} &=&\mathbf{i}\frac{\partial ^{2}L}{\partial z^{j}\partial z^{i}}%
dz^{i}\wedge dz^{j}+\mathbf{i}\frac{\partial ^{2}L}{\partial \overline{z}%
^{j}\partial z^{i}}dz^{i}\wedge d\overline{z}^{j}  \nonumber \\
&&+\mathbf{i}\frac{\partial ^{2}L}{\partial z^{j}\partial \overline{z}^{i}}%
dz^{j}\wedge d\overline{z}^{i}+\mathbf{i}\frac{\partial ^{2}L}{\partial
\overline{z}^{j}\partial \overline{z}^{i}}d\overline{z}^{j}\wedge d\overline{%
z}^{i}.  \label{3.6}
\end{eqnarray}

Let $\xi $ be the semispray given by (\ref{3.1}) and
\begin{eqnarray}
i_{\xi }\Phi _{L} &=&\mathbf{i}\xi ^{i}\frac{\partial ^{2}L}{\partial
z^{j}\partial z^{i}}dz^{j}-\mathbf{i}\xi ^{i}\frac{\partial ^{2}L}{\partial
z^{j}\partial z^{i}}\delta _{i}^{j}dz^{i}+\mathbf{i}\xi ^{i}\frac{\partial
^{2}L}{\partial \overline{z}^{j}\partial z^{i}}d\overline{z}^{j}-\mathbf{i}%
\overline{\xi }^{i}\frac{\partial ^{2}L}{\partial \overline{z}^{j}\partial
z^{i}}\delta _{i}^{j}dz^{i}  \nonumber \\
&&+\mathbf{i}\xi ^{i}\frac{\partial ^{2}L}{\partial z^{j}\partial \overline{z%
}^{i}}\delta _{i}^{j}d\overline{z}^{i}-\mathbf{i}\overline{\xi }^{i}\frac{%
\partial ^{2}L}{\partial z^{j}\partial \overline{z}^{i}}dz^{j}+\mathbf{i}%
\overline{\xi }^{i}\frac{\partial ^{2}L}{\partial \overline{z}^{j}\partial
\overline{z}^{i}}\delta _{i}^{j}d\overline{z}^{i}\mathbf{-}\mathbf{i}%
\overline{\xi }^{i}\frac{\partial ^{2}L}{\partial \overline{z}^{j}\partial
\overline{z}^{i}}d\overline{z}^{j}.  \label{3.7}
\end{eqnarray}
Since the closed K\"{a}hlerian form $\Phi _{L}$ on $TM$ is symplectic
structure, it is obtained
\begin{equation}
E_{L}=\mathbf{i}\xi ^{i}\frac{\partial L}{\partial z^{i}}-\mathbf{i}%
\overline{\xi }^{i}\frac{\partial L}{\partial \overline{z}^{i}}-L
\label{3.8}
\end{equation}
and hence
\begin{equation}
\begin{array}{ll}
dE_{L}+\wedge ^{a}\omega _{a}= & \mathbf{i}\xi ^{i}\frac{\partial ^{2}L}{%
\partial z^{j}\partial z^{i}}dz^{j}-\mathbf{i}\overline{\xi }^{i}\frac{%
\partial ^{2}L}{\partial z^{j}\partial \overline{z}^{i}}dz^{j}-\frac{%
\partial L}{\partial z^{j}}dz^{j} \\
& +\mathbf{i}\xi ^{i}\frac{\partial ^{2}L}{\partial \overline{z}^{j}\partial
z^{i}}d\overline{z}^{j}-\mathbf{i}\overline{\xi }^{i}\frac{\partial ^{2}L}{%
\partial \overline{z}^{j}\partial \overline{z}^{i}}d\overline{z}^{j}-\frac{%
\partial L}{\partial \overline{z}^{j}}d\overline{z}^{j}+\wedge ^{a}\omega
_{a}.
\end{array}
\label{3.9}
\end{equation}
With respect to (\ref{1.3}), if (\ref{3.7}) and (\ref{3.9}) is equalized, it
is calculated as follows:
\begin{equation}
\begin{array}{l}
-\mathbf{i}\xi ^{i}\frac{\partial ^{2}L}{\partial z^{j}\partial z^{i}}dz^{j}-%
\mathbf{i}\overline{\xi }^{i}\frac{\partial ^{2}L}{\partial \overline{z}%
^{j}\partial z^{i}}dz^{j}+\frac{\partial L}{\partial z^{j}}dz^{j} \\
+\mathbf{i}\xi ^{i}\frac{\partial ^{2}L}{\partial z^{j}\partial \overline{z}%
^{i}}d\overline{z}^{j}\mathbf{+}\mathbf{i}\overline{\xi }^{i}\frac{\partial
^{2}L}{\partial \overline{z}^{j}\partial \overline{z}^{i}}d\overline{z}^{j}+%
\frac{\partial L}{\partial \overline{z}^{j}}d\overline{z}^{j}=\wedge
^{a}\omega _{a}
\end{array}
\label{3.10}
\end{equation}
Now, let the curve $\alpha :\mathbf{C}\rightarrow TM$ be integral curve of $%
\xi ,$ which satisfies equations
\begin{equation}
\begin{array}{l}
-\mathbf{i}\left[ \xi ^{j}\frac{\partial ^{2}L}{\partial z^{j}\partial z^{i}}%
+\stackrel{.}{\xi }^{i}\frac{\partial ^{2}L}{\partial \stackrel{.}{z}%
^{j}\partial z^{i}}\right] dz^{j}+\frac{\partial L}{\partial z^{j}}dz^{j} \\
+\mathbf{i}\left[ \xi ^{j}\frac{\partial ^{2}L}{\partial z^{j}\partial
\stackrel{.}{z}^{i}}+\stackrel{.}{\xi }^{j}\frac{\partial ^{2}L}{\partial
\stackrel{.}{z}^{j}\partial \stackrel{.}{z}^{i}}\right] d\stackrel{.}{z}^{j}+%
\frac{\partial L}{\partial \stackrel{.}{z}^{j}}d\stackrel{.}{z}^{j}=\wedge
^{a}\omega _{a}
\end{array}
\label{3.11}
\end{equation}
where $\omega _{a}=(\omega _{a})_{j}$ $dz^{j}+(\stackrel{.}{\omega }_{a})_{j}%
\stackrel{.}{dz}^{j}$and the dots mean derivatives with respect to the time.
We infer the equations
\begin{equation}
\frac{\partial L}{\partial z^{i}}-\mathbf{i}\frac{\partial }{\partial t}%
\left( \frac{\partial L}{\partial z^{i}}\right) =\wedge ^{a}(\omega
_{a})_{i},\,\,\,\,\,\,\frac{\partial L}{\partial \stackrel{.}{z}^{i}}+%
\mathbf{i}\frac{\partial }{\partial t}\left( \frac{\partial L}{\partial
\stackrel{.}{z}^{i}}\right) =\wedge ^{a}(\stackrel{.}{\omega }_{a})_{i}.
\label{3.12}
\end{equation}
Thus, by \textit{complex Euler-Lagrange equations} \textit{with constraints }%
we may call the equations obtained in (\ref{3.12}) on K\"{a}hlerian manifold
$TM.$ Then the quartet $(TM,\Phi _{L},\xi ,\overline{\omega })$ is named
\textit{mechanical system with constraints}.\newline

\textbf{Conclusion}

Finally, considering the above, complex analogous of the geometrical and
mechanical meaning of constraints given in \cite{deleon2,kiehn} may be
explained as follows.

\textbf{1)} Let $\overline{\omega }$ be a system of constraints on
K\"{a}hlerian manifold $TM.$ Then it may be defined a distribution $D$ on $%
\overline{\omega }$ as follows.
\begin{equation}
D(x)=\{\left. \xi \in T_{x}TM\right| \,\omega _{a}(\xi )=0,\,\mbox{ for all }%
\;a,\,1\leq a\leq r\}  \label{5.1}
\end{equation}

Thus $D$ is $(2m-r)$ dimensional distribution on $TM.$ In this case, a
system of complex constraints $\overline{\omega }$ is called holonomic, if
the distribution $D$ is integrable; otherwise we call $\overline{\omega }$
anholonomic. \thinspace Hence, $\overline{\omega }$ is holonomic if and only
if the ideal $\rho $ of $\wedge TM$ generated by $\overline{\omega }$ is a
differential ideal$.$ Obviously (\ref{3.12}) holds for holonomic as well as
anholonomic constraints. For a system of holonomic constraints, the motion
lies on a specific leaf of the foliation defined by $D.$

\textbf{2) }From (\ref{1.3}) it is obtained equalities of
\begin{equation}
0=(i_{\xi }\Phi )(\xi )=dE_{L}(\xi )=\xi (E_{L}),  \label{5.2}
\end{equation}

Therefore, the Lagrangian energy $E_{L}$\ on K\"{a}hlerian manifold $TM$ for
a solution $\alpha (t)$ of (\ref{3.12}) is conserved.

\thinspace \thinspace \thinspace \thinspace \thinspace \thinspace \thinspace
\thinspace \thinspace \thinspace \thinspace \thinspace \thinspace \thinspace
\end{titlepage}
\end{document}